\documentclass[11pt]{paper}

\usepackage{amssymb,amsmath,enumerate,theorem,epsfig,rotating,graphics,changebar,eepic}%t1enc
\usepackage{amsfonts, booktabs}
\usepackage{a4,latexsym,parskip}
\usepackage{hyperref}

\hypersetup{pdfauthor=MSAS} \hypersetup{pdftitle=DS}

\newcommand{\Aff}{\mathbb{A}}
\newcommand{\PP}{\mathbb{P}}

\newcommand{\C} [1][]{\mathbb{C}^{#1}}
\newcommand{\Q} [1] []{\mathbb{Q}_{#1}}
\newcommand{\N} [1][] {\mathbb{N}_{#1}}

\newcommand{\F}{\mathbb{F}}

\newcommand{\OO}{\mathcal{O}}

\newcommand{\qed}{\hfill \ensuremath{\Box}}
\newcommand{\lcm}{\mathrm{lcm}}

\theoremstyle{break} \newtheorem{Theorem}{Theorem}%[section]
\newtheorem{Proposition}[Theorem]{Proposition}
\newtheorem{Lemma}[Theorem]{Lemma}

\newtheorem{Definition}[Theorem]{Definition}

\newtheorem{Remark}[Theorem]{Remark}

\newtheorem{Conjecture}[Theorem]{Conjecture}

\begin{document}

\setlength{\unitlength}{1cm}

\title{Unirational Surfaces on the Noether Line}
\author{Christian Liedtke, Matthias Sch\"utt}
\date{\today}

\maketitle

\abstract{
We show that among simply connected surfaces of general type unirationality
is a common feature, even when fixing the positive characteristic or numerical invariants.
To do so, we construct unirational Horikawa surfaces in abundance.}

\keywords{unirationality, L\"uroth problem, Noether--Horikawa surface}

\textbf{MSC(2000):} 14J29, 14M20

\section{Introduction}

The L\"uroth problem asks whether a variety that is dominated by a rational variety is rational itself, i.e., whether unirationality implies rationality.
This is known to be true for curves (L\"uroth) and for complex surfaces (Castelnuovo).
However, it is also known to be wrong in dimension $\geq3$ thanks to 
counterexamples by Fano--Iskovskikh--Manin \cite{IM} and Clemens--Griffiths \cite{CG}.
Zariski \cite{Z} constructed unirational but non-rational surfaces over
fields of positive characteristic.
Artin and Mumford \cite{AM} gave unirational examples of dimension $\geq 3$ in characteristic $\neq 2$ that cannot be rational if the characteristic is zero or the dimension is 3.

On the other hand, one can ask what kind of varieties can be unirational:
over the complex numbers a unirational variety is always of Kodaira dimension $\kappa=-\infty$.
Moreover, Serre \cite{se}
has shown that the fundamental group is finite in general, and even 
trivial in characteristic zero.

This leads to the question among what kind of surfaces in positive characteristic, say with trivial fundamental group, we can find unirational ones.
Clearly, the most ''non-rational'' class of surfaces are those of general type, i.e., of Kodaira dimension $\kappa=2$.

Among them we pick a very prominent series, the so-called
{\em Horikawa surfaces}.
We recall that minimal surfaces of general type
fulfill Noether's inequality 
$$
K^2\,\geq\,2p_g\,-\,4\,.
$$
By definition, Horikawa surfaces are those for which equality holds,
which is why they are also known as
{\em surfaces on the Noether line}.
These surfaces have a particularly simple structure: their canonical map
exhibits their canonical models as double covers
of minimal rational surfaces,
cf.~Theorem~\ref{thm:Hori} and \cite[Section VII.9]{bhpv}.

These double covers are known to be simply connected.
In characteristic $2$ it can happen that the canonical map becomes
inseparable. Then the surfaces are necessarily unirational,
cf. \cite[Section 5]{lie}.
However, since the canonical map of a Horikawa surface is always of degree $2$, 
this inseparability phenomenon can happen in characteristic $2$ only.

Over fields of positive characteristic $p\geq3$, we therefore face the 
following problems:

\begin{enumerate}
\item Do there exist unirational Horikawa surfaces in characteristic $p\geq3$?
\item If so, is the characteristic bounded?
\item If not, is the number of possible characteristics 
  for unirational Horikawa surfaces with fixed $p_g$ finite?
\item Is the number of possible geometric genera $p_g$
  for unirational Horikawa surfaces with fixed characteristic finite?
\end{enumerate}

Shioda showed that unirational surfaces are supersingular, 
i.e., that their Picard numbers are equal to their second Betti numbers -
even if $h^{0,2}$ or $h^{2,0}$ are non-zero.
Moreover, he also conjectured the converse for simply connected surfaces, 
cf. \cite{sh2} and Conjecture \ref{conj} below.
So far, the Picard numbers of Horikawa surfaces have only been studied in the 
complex case by Persson \cite{Per}, cf. also Section~\ref{s:class}.
Hence these questions can also be viewed as a characteristic $p$ extension
of Persson's work.

As to the first and the second question, we prove existence of unirational 
Horikawa surfaces in $99,999985\%$ positive characteristics.
More precisely, we show

\begin{trivlist}
\item[]{\bf Theorem \ref{rem:dense}}\\
\emph{For all primes $p$  outside a set $\mathcal{P}$ of density $5\cdot2^{-25}$, 
there do exist unirational Horikawa surfaces in characteristic
$p$.}
\end{trivlist}

Hence the second question has a negative answer, but even the refined third
question has a negative answer:

\begin{trivlist}
\item[]{\bf Theorem \ref{horikawa}}\\
\emph{Let $g$ be an even integer.
Then there exists an arithmetic progression $P_g$ of primes
such that for every $p\in P_g$ there exists a unirational
Horikawa surface with $p_g=g$ in characteristic $p$.}
\end{trivlist}

Using Delsarte surfaces, we obtain much more detailed results for
Horikawa surfaces with $p_g=3$ and $p_g=6$, cf. 
Theorem \ref{thm:p_g=3} and Remark \ref{rem:p_g=6}.

The fourth question asked what happens when fixing
the characteristic:

\begin{trivlist}
\item[]{\bf Theorem \ref{question4}}\\
\emph{Let $p\geq3$ be a prime with $p\not\equiv1\mod15$.
Then there exists an infinite set $G_p$ of positive
integers such that there exist unirational Horikawa
surfaces in characteristic $p$
with $p_g=g$ for all $g\in G_p$. }
\end{trivlist}

These results show that unirational Horikawa surfaces exist
in abundance - even when fixing the characteristic $p$ or the geometric
genus $p_g$.
In Section \ref{s:singular} we show moreover that the unirationality
cannot be explained by singular fibrations arising from pluricanonical
systems.

As an interesting byproduct we construct a five-dimensional family
of supersingular Horikawa surfaces with $p_g=3$ in characteristic $p=7$
where all surfaces have Artin invariants $3\leq\sigma\leq4$.
Unlike for supersingular K3 surfaces, this implies that the moduli space of supersingular Horikawa surfaces with $p_g=3$ in characteristic $7$ does not admit a stratification with moduli dimension equalling Artin invariant minus one.

%In this paper, we give an affirmative answer to the first question for a great density of %characteristics (see Rem.~\ref{rem:dense}). The other two questions have a negative answer %for all even $p_g$ and $p_g=3$ by Thm.s \ref{thm:p_g=3} and \ref{horikawa}.

The paper is organised as follows: The next section recalls Enriques' and Horikawa's classification and the basic invariants and notions. Section \ref{s:fam} introduces some families of unirational Horikawa surfaces that arise naturally. In sections \ref{s:p_g=3} and \ref{s:even}, we construct the surfaces to answer questions 1-4. The final section discusses whether unirationality can be explained through singular fibrations arising from pluricanonical systems.

\section{Enriques' and Horikawa's classification}
\label{s:class}

In this section, we recall the main results about Horikawa surfaces and in particular Enriques' and Horikawa's classification.

\begin{Definition}
A minimal surface of general type is called \emph{Horikawa surface} if it lies on the Noether line
\[
K^2=2\,p_g-4.
\]
\end{Definition}

The maybe most classical example for an Horikawa surface is a double octic: the double cover of $\PP^2$ branched along an octic curve $C$. Here we can allow the curve $C$ to have isolated ade-singularities. Then we have to consider the minimal resolution of the resulting ADE-singularities on the double cover. 

The same ideas apply to the construction of the other Horikawa surfaces. 
Over the complex numbers, they were classified by Enriques and Horikawa \cite{Hori}. 
The techniques were carried over to positive characteristic by the first author in \cite{lie}. 
We denote by $\F_d$ the $d$-th Hirzebruch surface.

\begin{Theorem}[Enriques, Horikawa]
\label{thm:Hori} 
Let $X$ be a Horikawa surface. Via the canonical map $\phi_1$, $X$ is the double cover of a rational surface $S:=\phi_1(X)$ in $\PP^{p_g-1}$ (possibly singular).
The  following cases occur:
 \begin{enumerate}
 \item If $S$ is a smooth surface, then we have the
   following possibilities
   \begin{enumerate}[(i)]
   \item  $S\cong\PP^2$ and $p_g=3$ (double octic).
   \item $S\cong\PP^2$ and $p_g=6$ (the branch locus is a curve of degree 10).
   \item $S\cong\F_d$ and $p_g\geq\max\{d+4, 2d-2\}$
     and $p_g-d$ is even.
   \end{enumerate}
 \item If $S$ is not smooth, then it is the cone over the rational normal
   curve of degree $d:=p_g-2$.
   The minimal desingularisation of $S$ is the Hirzebruch
   surface $\F_{d}$ and
   $4\leq p_g\leq 6$.
 \end{enumerate}
\end{Theorem}

It follows from this description that complex Horikawa surfaces are 
topologically simply connected \cite[Thm.~3.4]{Hori}
and that they are algebraically simply connected 
with reduced Picard schemes in positive 
characteristic \cite[Prop.~3.7]{lie}. 

Given the geometric genus $p_g$ of a complex Horikawa surface, we compute the remaining Hodge number 
\[
h^{1,1}=8\,p_g+14.
\]
In \cite{Per}, Persson investigated complex Horikawa surfaces with maximal Picard number $\rho=h^{1,1}$. He proved existence whenever $p_g\not\equiv 5\mod 6$ or $p_g=5$.

In positive characteristic, however, Igusa showed that the Picard number can be as big as $b_2$. Such surfaces are often called \emph{supersingular}. In \cite{sh2}, Shioda formulated the following

\begin{Conjecture}\label{conj}
Let $k$ an algebraically closed field of positive characteristic. Let $X$ be an algebraically simply connected surface over $k$. Then $X$ is unirational if and only if it is supersingular. 
\end{Conjecture}

The "only if"-implication of the conjecture holds almost trivially true, cf. \cite{sh2}.
However, in practice it is often easier to verify that a given surface is supersingular than unirational. For instance, one can start by computing the characteristic polynomial of Frob$_p$ on the second \'etale cohomology and apply the Tate conjecture \cite{tate}. We pursued this approach without success for the reductions of several families of complex Horikawa surfaces of maximal Picard number borrowed from Persson's construction. In the next sections, we thus present surfaces where we can prove unirationality in a direct and systematic way.

As a side remark, we note a very different behaviour for K3 surfaces. We regard them close to Horikawa surfaces as they often can be realised as double sextics. Complex K3 surfaces of maximal Picard number 20 are called singular. They are necessarily defined over a number field. The supersingularity of the reductions does only depend on the splitting behaviour in an imaginary-quadratic field. The interested reader is referred to Shioda-Inose's paper \cite{SI} for details.

\section{Families of unirational Horikawa surfaces}
\label{s:fam}

In this section, we give an ad-hoc example of a family of unirational Horikawa surfaces with $p_g=3$ in characteristic $7$. It mimics a construction of Pho and Shimada for unirational K3 surfaces in characteristic $5$ \cite{ps}.

The Horikawa surfaces are realised as double octics over an algebraically closed field $k$ of characteristic $7$. The branch curve $C$ is affinely given in $\mathbb{A}^2$ by
\begin{eqnarray}\label{eq:C}
C:\;\;\; y^7 = f(x).
\end{eqnarray}
Here $f\in k[x]$ is a polynomial of degree $8$ such that its formal derivative has no multiple roots in $k$. In consequence, $C$ is smooth outside seven $a_6$ singularities where the formal derivative of $f$ vanishes. 

Let $Y$ be the double cover of $\PP^2$ branched along $C$. Denote the minimal resolution of the seven $A_6$ singularities of $Y$ by $X$. By construction, $X$ is a unirational Horikawa surface with $p_g=3$ for any such polynomial $f$. As we can rescale and apply M\"obius transformations, we obtain a five-dimensional family of unirational Horikawa surfaces.

Recall that every Horikawa surface with $p_g=3$ arises as a double cover of $\PP^2$ branched
along an octic curve with at worst ade-singularities. Hence the moduli space of all
such surfaces is of dimension $36$.

\begin{Proposition}
Let $X$ be a Horikawa surface as above. Then for some $3\leq\sigma\leq 4$
\[
\text{disc } NS(X)=-7^{2\sigma}.
\]
\end{Proposition}

\emph{Proof:}
The existence of such an invariant $\sigma\geq 1$ was worked out in \cite{lie2} and \cite{SS}. 
It is not difficult to see that the minimal desingularisation $X$ of $Y$ for every
$Y$ in our family lifts over the Witt ring.
In particular, the Fr\"olicher spectral sequence degenerates at $E_1$ and the Hodge numbers 
in characteristic $p=7$ and zero coincide.
Hence $h^0(\Omega_X^1)=0$ and since ${\rm Pic}(X)$ is reduced, it follows from
\cite[Prop.~II.5.16]{ill} that the second crystalline cohomology of $X$ is torsion free.
From \cite[Rem.~II.5.21]{ill} or \cite[Prop.~6.1]{lie2} we obtain $\sigma\geq p_g=3$.
We will prove $\sigma\leq 4$ by exhibiting a suitable sublattice of finite index in $NS(X)$.

The exceptional curves of the minimal resolution $X\to Y$ generate a sublattice 
\[
N = 7 A_6(-1) \subset NS(X).
\]
As $b_2(X)=10\,p_g(X)+14=44$, we shall exhibit another sublattice $L$ of $NS(X)$ of rank two which is orthogonal to $N$. For this, consider the line $\ell$ at infinity with respect to the chart of $C$ in (\ref{eq:C}). Then $\ell$ meets $C$ in a single point of multiplicity 8. Hence the pull-back of $\ell$ under the composition $X \to Y \to \PP^2$ splits into two curves
\[
\pi^* \ell = \ell_1 +\ell_2
\]
which meet in a single point with multiplicity 4. Since $(\pi^* \ell)^2=2$, it follows that $\ell_i^2=-3$. Hence the curves $\ell_i$ span an indefinite sublattice $L$ of $NS(X)$ with intersection form
\[
\begin{pmatrix} -3 & 4\\4 & -3\end{pmatrix}.
\]
As $C$ is smooth along $\ell$, $L$ is orthogonal to $N$. Hence $N+L$ is a sublattice of $NS(X)$ of finite index and discriminant $-7^8$. Thus $\sigma\leq 4$. \qed

\begin{Remark}
As we have exhibited a five-dimensional family of (unirational) Horikawa surfaces with these discriminants, we deduce that
the moduli space of supersingular Horikawa surfaces with $p_g=3$ in characteristic $7$ does not admit a stratification with moduli dimension equalling invariant $\sigma$ minus one.

This result is contrary to the situation for supersingular K3 surfaces in arbitrary characteristic as investigated by Artin  \cite{A}. For supersingular K3 surfaces, $\sigma$ is called Artin invariant and ranges from 1 to 10, giving the moduli dimension $\sigma-1$.
\end{Remark}

The analogous construction can be applied to double covers of $\PP^2$ branched along a curve of degree 10. This gives rise to a seven-dimensional family of unirational Horikawa surfaces of geometric genus $p_g=6$ in characteristic $3$.

\section{Unirational Horikawa surfaces with $p_g=3$}
\label{s:p_g=3}

In this section we shall only construct single Horikawa surfaces over $\Q$. However, we will show that their reductions are unirational for all primes in some arithmetic progression. This gives a negative answer to the second question stated in the introduction. The surfaces will have geometric genus $p_g=3$, but the same approach works also for $p_g=6$. 

The main idea is as follows: We will exhibit all those Horikawa surfaces which are at the same time double octics and Delsarte surfaces. After Shioda \cite{dels}, we shall then use that Delsarte surfaces are dominated by Fermat surfaces. Thanks to work by Shioda and Katsura \cite{SK}, we gain complete knowledge about unirationality.

\begin{Definition}\label{Def}
A surface $X\subset\PP^3$ is called \emph{Delsarte surface} if it is defined by a homogeneous polynomial of four monomials
\[
X: \;\;\; \sum_{i=0}^3 x_0^{a_{i0}}\,x_1^{a_{i1}}\,x_2^{a_{i2}}\,x_3^{a_{i3}} = 0
\]
satisfying the following conditions:
\begin{enumerate}[(i)]
\item for each $j$, there is some $i$ with $a_{ij}=0$;
\item $\det\, (a_{ij})_{i,j} \neq 0$.
\end{enumerate}
\end{Definition}

Shioda considered Delsarte surfaces because they are always dominated by some Fermat surface
\[
S = \{y_0^m+y_1^m+y_2^m+y_3^m=0\}\subset\PP^3.
\]
Here the appropriate degree $m$ of the Fermat surface can be computed purely in terms of the matrix $A=(a_{ij})_{i,j}$. Let $A^*=(a_{i,j}^*)_{i,j}$ denote the cofactor matrix of $A$, i.e., $A\cdot A^* = (\det A)\cdot \mathbf{1}$. Set $\delta=\gcd(a_{i,j}^*)$ and $d=(\det A)/\delta$. Then the Delsarte surface $X$ given by $A$ is dominated by the Fermat surface of degree $d$. Writing $b_{i,j}=a_{i,j}^*/\delta$, the covering map is given as
\begin{eqnarray*}
\phi:\;\;\;\;S\;\, & \to & \;\;\;\;\;X\\
(y_i) & \mapsto & (x_i=\prod_{j=0}^3 y_j^{b_{i,j}}).
\end{eqnarray*}
Sometimes a smaller degree suffices for the covering Fermat surface.  This issue will be important for Thm.~\ref{rem:dense}. We will investigate the question of the precise degree for one example in the sequel of Thm.~\ref{thm:p_g=3}. A general statement can be found in \cite[Def.~1 \& Lem.~2]{Crelle}.

Fermat surfaces are very well-understood. The study of their arithmetic goes back to Weil who computed the $\zeta$-function over finite fields in terms of Jacobi sums. In \cite{SK}, Shioda and Katsura proved Conjecture \ref{conj} for Fermat surfaces:

\begin{Theorem}[Shioda--Katsura]
The Fermat surface of degree $m$ is unirational in characteristic $p$ if and only if it is supersingular. Equivalently, there is some $\nu\in\N$ with
\begin{eqnarray}\label{eq:nu}
p^\nu\equiv -1 \mod m.
\end{eqnarray}
\end{Theorem}

One prominent feature of Fermat surfaces is the motivic decomposition of their cohomology in terms of some character group. This decomposition allows the explicit calculation of Picard number $\rho(S)$ and Lefschetz number
\[
\lambda(S)=b_2(S)-\rho(S)
\]
over any algebraically closed field of characteristic coprime to the degree $m$. 

A Delsarte surface arises from the dominating Fermat surface as the quotient by some group action. Hence, the transcendental part of the cohomology of the quotient surface can be identified with the transcendental part of the original surface which is invariant under the group. Thus we can calculate the Lefschetz number of any Delsarte surface and consequently also the Picard number. Clearly, the Delsarte surface is unirational if the Fermat surface is.

With this knowledge, we find unirational Horikawa surfaces as double octics which at the same time are Delsarte surfaces:

\begin{Theorem}
\label{thm:p_g=3}
Let $X$ be the minimal resolution of the double cover of $\PP^2$ branched along one of the following octic curves $C$. Then $X$ is dominated by the Fermat surface of degree $m$ as in the table. In particular, $X$ is unirational in characteristic $p$ if and only if (\ref{eq:nu}) holds for the respective degree $m$:

\pagebreak

$$
\begin{array}{|c|c|}
\toprule 
C & m\\
\midrule
\midrule
x^8+y^8+z^8=0 & 8\\
\midrule
x^7\,y+y^7\,z+z^7\,x=0 & 86\\
\midrule
x\,(x^7+y^7+z^7)=0 & 14\\
\midrule
x\,(x^6\,y+y^6\,z+z^6\,x)=0 & 62\\
\midrule
x\,y\,(x^6+y^6+z^6)=0 & 12\\
\midrule
x\,y\,(x^5\,y+y^5\,z+z^5\,x)=0 & 42\\
\midrule
x\,y\,z\,(x^5+y^5+z^5)=0 & 10\\
\midrule
x\,y\,z\,(x^4\,y+y^4\,z+z^4\,x)=0 & 26\\
\midrule
x^8+z\,y^7+x\,y\,z^6=0 & 82\\
\midrule
x\,(x^7+z\,y^6+x\,y\,z^5)=0 & 58\\
\midrule
z\,(x^6\,z+y^5\,x^2+z^5\,y^2)=0 & 44\\
\midrule
x^7\,z+y^7\,x+z^6\,y^2 = 0 & 74\\
\midrule
x\,z\,(x^5\,z+y^5\,x+z^4\,y^2) = 0 & 34\\
\bottomrule
\end{array}
$$
\end{Theorem}

Over $\C$, all curves in Thm.~\ref{thm:p_g=3} have only isolated ade-singularities. Hence we obtain Horikawa surfaces. In positive characteristic, there are only finitely many exceptions where the singularities of a given curve degenerate. In the cases at hand, $C$ is dominated by the Fermat curve of degree $m$ (cf.~the example below). Since this Fermat curve has good reduction outside the primes dividing $m$, the same holds for $C$. %The primes dividing $m$ 
These primes are excluded by the condition (\ref{eq:nu}).

For one Horikawa surface in Thm.~\ref{thm:p_g=3}, we will now sketch how to obtain the given degree of the covering Fermat surface. Consider the  second  octic curve
\[
C=\{x^7\,y+y^7\,z+z^7\,x=0\}\subset\PP^2.
\]
We start by determining the covering Fermat curve along the  lines of Def.~\ref{Def}. Let $A$ be the matrix of exponents. Then $A$ has determinant $344$. Its cofactor matrix is
\[
A^* = \begin{pmatrix} 49 & -7 & 1\\1 & 49 & -7\\-7 & 1 & 49\end{pmatrix},
\]
so $C$ is covered by the Fermat curve of degree $344$. Since we work in projective space, the rational map $\phi$ stays the same when we add constants to the columns of $A^*$. After adding $7$ to all columns, the coefficients become divisible by $8$. Division yields the matrix
\[
B = \begin{pmatrix} 7 & 0 & 1\\1 & 7 & 0\\0 & 1 & 7\end{pmatrix}.
\]
The corresponding rational map $\phi$ takes the Fermat curve of degree $43=344/8$ to $C$. One can easily complement $\phi$ to a rational map from a Fermat surface to the Horikawa surface associated to $C$.
Since the Horikawa surface is exhibited as a double cover, the dominating Fermat surface has degree $86$. Expressed through the matrix of exponents, the covering map corresponds to
\[
\begin{pmatrix}
43 & 7 & 7 & 7\\
0 & 14 & 0 & 2\\
0 & 2 & 14 & 0\\
0 & 0 & 2 & 14
\end{pmatrix}.
\]

\begin{Remark}
There are further octic curves with only ade-singularities which are defined by a three-term polynomial. They give rise to more unirational Horikawa surfaces as double cover. However, the degrees of the dominating Fermat surface are not essentially new in the following sense: The degree is always a multiple of some degree in the previous table or of $m=6$ which will turn up as Fermat degree in the next section (cf.~Thm.~\ref{horikawa}). The precise  degrees occurring will be relevant for Thm.~\ref{rem:dense}.
\end{Remark}

\begin{Remark}
\label{rem:p_g=6}
The same construction applies to Horikawa surfaces as double covers of $\PP^2$ branched along a curve of degree $10$. Allowing ade-singularities, we find many further unirational Horikawa surfaces. There is one essentially new degree for the dominating Fermat surfaces: $m=146$. It is achieved by the following curve of degree 10:
\[
\{x^9\,y+y^9\,z+z^9\,x=0\}\subset\PP^2.
\]
\end{Remark}

\section{Unbounded invariants}
\label{s:even}

In the previous section we have seen unirational Horikawa surfaces 
with $p_g=3$ for a large density 
of characteristics.
We now construct unirational Horikawa surfaces 
in arithmetic progressions of characteristics for every
even value of $p_g$.

\begin{Lemma}
\label{shiodalemma}
Let $a,b$ be positive integers and set $m:=\lcm(a,b)$.
Then the double cover of $\Aff^2$
$$
  X_{a,b}\,:\,z^2\,=\,(x^a-1)\cdot(y^b-1)
$$
defines a unirational surface in every characteristic
$p\geq3$ for which there exists a $\nu\in\N$ 
with $p^\nu\equiv-1\mod m$.
\end{Lemma}

\emph{Proof:}
For every integer $n\geq3$ we consider the hyperelliptic curve
$$
C_n\,:\,v^2\,=\,1\,-\,u^n
$$
together with its involution $\imath:(u,v)\mapsto(u,-v)$.
For $m=\lcm(a,b)$ the product $C_m\times C_m$ dominates $C_a\times C_b$.
By abuse of notation we also denote by $\imath$ the involution on
these products given by $(u_1,v_1,u_2,v_2)\mapsto(u_1,-v_1,u_2,-v_2)$.
Then, $Y_m:=(C_m\times C_m)/\imath$ dominates 
$(C_a\times C_b)/\imath$
and  $Y_m$ is unirational in characteristic $p$ if and only if there exists
a $\nu$ such that $p^\nu\equiv-1\mod m$ by \cite[Lemma 1.2]{sh2}.
In particular, $(C_a\times C_b)/\imath$ is unirational in these
characteristics. An
elementary calculation shows that this surface is birational to
$X_{a,b}$.
\qed

\begin{Theorem}
\label{horikawa}
Let $g$ be an even integer.
Then there exists an arithmetic progression $P_g$ of primes
such that for every $p\in P_g$ there exists a unirational
Horikawa surface with $p_g=g$ in characteristic $p$.
\end{Theorem}

\emph{Proof:}
We define
$$
\begin{array}{ccl}
  M_g &:=&\{\,\lcm(5,g+1),\, \lcm(5,g+2),\,
  \lcm(6,g+1),\,\lcm(6,g+2)\,\},\\
  \rule[-3mm]{0mm}{8mm}
  P_g &:=&\{\,p\in\N\,|\,p\geq3 \mbox{ a prime s.t. }
  \exists\nu\in\N, m\in M_g\,:\,p^\nu\equiv-1\mod m\,\},
\end{array}
$$
and note that $P_g$ is the 
union of primes in arithmetic progressions.

We use the surfaces $X_{a,b}$ from Lemma \ref{shiodalemma}.
Let $a'$ (resp.~$b'$) be the smallest even integer bigger 
than or equal to $a$ (resp.~$b$).
We think of $\Aff^1\times\Aff^1$ as contained in 
$\PP^1\times\PP^1$.
If we consider $(x^a-1)(y^b-1)$ as a global section of 
$\OO_{\PP^1\times\PP^1}(a',b')$ then its zero set 
is the union of $a'+b'$ lines intersecting transversally.
Hence a double cover branched over this divisor is a surface
with only Du~Val singularities of Type $A_1$.

In case $a'=6$ (i.e., $a=5$ or $a=6$) and 
$b'=g+2$ (i.e., $b=g+1$ or $b=g+2$) this surface
is the canonical model of a Horikawa surface 
with $p_g=g$, cf. \cite[Sect.~VII.9]{bhpv} or
\cite[Thm.~3.3]{lie}.
By Lemma \ref{shiodalemma} this surface is unirational in
every characteristic $p$ for which there exists a
$\nu$ with $p^\nu\equiv-1\mod\lcm(a,b)$.
Hence for every $p\in P_g$ we have constructed a
unirational Horikawa surface with $p_g=g$ in characteristic $p$.
\qed

\begin{Theorem}\label{question4}
Let $p\geq3$ be a prime with $p\not\equiv1\mod15$.
Then there exists an infinite set $G_p$ of positive
integers such that there exist unirational Horikawa
surfaces in characteristic $p$
with $p_g=g$ for all $g\in G_p$. 
\end{Theorem}

\emph{Proof:}
If $p$ is a prime as in the statement then there exists
an infinite number of $\nu$ such that
$p^\nu\equiv-1\mod5$ or $p^\nu\equiv-1\mod6$.
For every such $\nu$, we obtain a unirational
Horikawa surface with $p_g=(p^\nu+1)-2$ in characteristic
$p$ using the examples constructed in the proof of 
Theorem \ref{horikawa}.
\qed

\begin{Theorem}\label{rem:dense}
For all primes $p$ outside a set $\mathcal{P}$ of density $5\cdot2^{-25}$, 
there do exist unirational Horikawa surfaces in characteristic
$p$.
\end{Theorem}

\emph{Proof:}
Consider all the surfaces constructed in this and the preceding section. Their unirationality mod $p$ does only depend on the respective Fermat degree $m$. Condition (\ref{eq:nu}) does not distinguish between odd modulus $m$ and $2m$, so we list odd moduli whenever applicable. 
%Here we can neglect multiples and consider $p$ instead of $2p$. Hence we have to take the following degrees into account:
We constructed unirational Horikawa surfaces for the following moduli:
\[
\mathcal{M}=\{3,5,7,8,13,17,29,31,37,41,43,44,73\}.
\]
Hence the set $\mathcal{P}$ can be defined as follows:
\[
\mathcal{P} = \{p; \; \forall \,m\in\mathcal{M}, \nu\in\N:\;\; p^\nu\not\equiv -1\mod m\}.
\]
Working out the congruence conditions for each modulus $m\in\mathcal{M}$, one easily computes the claimed density. \qed

\begin{Remark}
We have thus constructed unirational Horikawa surfaces for a density of $99.999985\%$ of characteristics.
The smallest characteristic for which we have no example is $p=67665781$.
\end{Remark}

\begin{Remark}
For all surfaces in this and the previous section, we can compute the zeta function by motivic decomposition. For the dominating Fermat varieties, this decomposition is due to Weil \cite{W}. 
Then one computes the submotive invariant under the automorphism group corresponding to the covering map to obtain the essential factor of the zeta function.
% through the dominating Fermat surface via motivic decomposition and group invariance.
\end{Remark}

\section{A remark on singular fibrations}
\label{s:singular}

An interesting question is what \emph{causes} unirationality.
Clearly, unirational surfaces are covered by possibly singular rational curves.
On the other hand, if a surface admits a fibration onto a curve whose
generic fibre is a possibly singular geometrically rational curve then the surface
in question is uniruled.
In case the base of this fibration is $\PP^1$ 
(in our case this would be automatic since our surfaces have $b_1=0$), 
the surface is even unirational.

For surfaces of general type, it is natural to investigate their 
pluricanonical systems.
Maybe the unirationality of Horikawa surfaces, say with fixed $p_g$,
can be explained by such singular fibrations arising from pluricanonical
systems?

\begin{Proposition}
\label{prop:singularfibration}
Let $X$ be a Horikawa surface in characteristic $p$ and $n\in\N$. 
Assume that some subsystem of $|nK_X|$ defines
a (possibly rational) map with $1$-dimensional image whose
generic fibre is a singular and geometrically rational curve.
Then
$$
 p\,\leq\,2\,+\,(n^2+n)\,(p_g(X)-2)
$$
holds true.
\end{Proposition}

\emph{Proof:}
Let $f:X\dashrightarrow C$ be a rational map with $1$-dimensional image
defined by some subsystem of $|nK_X|$.
In particular, if $F$ denotes the generic fibre of this fibration then
$A:=nK_X-F$ is effective.
Since $K_X$ is nef and $A$ is effective, we have $K_X F\leq nK_X^2$.
Moreover, we compute
$F^2=(nK_X-A)^2=n^2K_X^2-A(2nK_X-A)=n^2K_X^2-A(nK_X+F)\leq n^2K_X^2$
using numerical connectedness of pluricanonical divisors.
Using the adjunction formula and the definition of a 
Horikawa surface we get
$$
p_a(F)\,=\,1+\frac{1}{2}(F^2+K_XF)\,\leq\,
1\,+\,\frac{1}{2}(n^2+n)K_X^2\,=\,1\,+\,\frac{1}{2}(n^2+n)\,(2p_g(X)-4).
$$
We assumed the generic fibre of $f$ to be a singular curve, i.e., that
this curve is regular but not smooth over the function field of $C$.
By Tate's theorem on genus change in inseparable field extensions
$p-1$ divides $p_a(L)$, since we assumed $L$ to be geometrically rational.
In particular, $p_a(L)$ must be larger than or equal to $p-1$.
\qed

However, we have seen that even for fixed $p_g$ there may be
an infinite number of characteristics in which there exist
unirational Horikawa surfaces with this given $p_g$.
Proposition \ref{prop:singularfibration} tells us that
(possibly rational) fibrations with geometrically rational fibres
arising from $|nK_X|$ with $n$ bounded can only exist
in a finite number of characteristics.
Hence the unirationality of Horikawa surfaces in arbitrary
large characteristics is not related to the existence of
such singular fibrations arising
from a finite number of pluricanonical systems.
\vspace{0.5cm}

\textbf{Acknowledgement:} 
Much of this work was done while the second author was visiting the Universit\"at D\"usseldorf. We thank the Mathematisches Institut and the DFG Forschergruppe 
\emph{Classification of Algebraic Surfaces and Compact Complex Manifolds}
for the kind hospitality and support. The second author gratefully acknowledges funding from DFG under research grant Schu 2266/2-2. We thank the referee for helpful comments.

%\vspace{0.5cm}

%\pagebreak

\vspace{0.8cm}

\begin{center}
\begin{tabular}{llll}
Christian Liedtke &\hspace{1cm} && Matthias Sch\"utt\\
Mathematisches Institut&&& Department of Mathematical Sciences\\
Universit\"at D\"usseldorf&&&University of Copenhagen\\
Universit\"atsstra\ss e 1&&& Universitetspark 5\\
40225 D\"usseldorf&&& 2100 Copenhagen\\
Germany&&& Denmark\\
{\tt liedtke@math.uni-duesseldorf.de}&&&{\tt mschuett@math.ku.dk} 
\end{tabular}
\end{center}

\end{document}